\documentclass[3p,times]{elsarticle}

\usepackage{ecrc}
\usepackage{natbib}
\usepackage{latexsym,amsmath}
\usepackage{color}
\usepackage{pstricks,pst-slpe}

\volume{00}

\firstpage{1}

\journalname{Icarus}

\runauth{E. Verheylewegen et al.}

\jid{}

\usepackage{amssymb}

 \biboptions{comma,authoryear}

\usepackage[figuresright]{rotating}

\begin{document}

\begin{frontmatter}




\title{Coupled orbital-thermal evolution of Miranda}


\author[unamur]{E. Verheylewegen}
\author[obs]{\"O. Karatekin} 
\author[unamur]{B. Noyelles}

\address[unamur]{Namur Centre for Complex Systems, naXys, University of Namur, Namur, Belgium}
\address[obs]{Royal Observatory of Belgium, Brussels, Belgium}

\begin{abstract}
Miranda has a unusually high inclination ($I=4.338^\circ$), and its surface reveals signs of past endogenic activity. Investigations of the 
dynamical aspects of its orbital evolution suggest probable resonant processes, in particular with Umbriel, as an explanation for the present 
high inclination of Miranda. The tidal heating induced by gravitational interactions can lead to the rise of eccentricities and, consequently, 
to the increased dissipation of energy inside the satellite and higher internal temperatures. We study here the possible increase in eccentricities caused 
by orbital resonances and the resulting endogenic heating on Miranda taking into account its temperature dependent rheology. The coupled orbital-thermal evolution model was 
run with different rheological models and the thermal parameters starting form a cold thermal state, in radiative equilibrium with the environment. For the nominal parameters of 
the evolution scenarios studied, the resonances were not sufficient to rise neither the eccentricities nor the internal temperatures significantly. Lowest dissipation function $Q$
of around 100 and final eccentricity of $e\approx0.02$ were obtained during the resonance 3:1 with Umbriel.
\end{abstract}

\begin{keyword}
Uranus, satellites; Satellites, dynamics; Thermal histories; Resonances, orbital; moons, interior



\end{keyword}

\end{frontmatter}


\section{Introduction}
\label{Introduction}

The flyby of the giant planets and their moons by {\it Voyager 2} resulted in the most valuable planetary data set on the outer Solar System. For the Uranian System, the southern hemispheres
of the moons were the only enlightened part \citep{smith1986} but these images show that the surfaces of the main Uranian satellites present signs of both endogenic activity 
and of the impact environment in the early stages of the evolution \citep{brown1991}. The spacecraft {\it Voyager} made closest encounter flyby with Miranda, and was able to capture details
of the tectonic structures on the surface with relatively high resolution. The data from Miranda as well as Ariel show signs of endogenic resurfacing associated with cryo-volcanism process (\citealt{plescia1987,plescia1988}).

The moon Miranda is enigmatic because Miranda has a quite small size \citep{thomas1988} in comparison with other satellites of Uranus with complex  geological features  suggesting  potentially interesting geological history. The surface is composed of two types of fields: older craterised regions and regions called coronae \citep{strobell1987} showing signs of diapirism phenomenon \citep{pappalardo1997}. Following \citet{brown1991}  the thermal history of Miranda is divided into at least 2 distinct periods where the coronae structures  and diapirism  appeared in the last period.    
A probable explanation of the coronae structures is given by a tidal heating induced by gravitational interaction between satellites. The sine-qua none condition is the pumping  of eccentricities by resonance processes. For Miranda, \citet{dermott1988} estimate an increase of $20$ K with a pumped eccentricity of $0.1$. \citet{tittemore1990} evaluate the tidal heating of Miranda induced by the 3:1 mean-motion resonance with Umbriel. They observe large variations in eccentricities during chaotic stages of the evolution
but the final eccentricity of Miranda is not sufficiently maintained if the tidal heating is the only considered process. Another important dynamical element of the Uranian System pointed out by \citet{tittemore1989} is the current high inclination of Miranda ($I=4.338^\circ$) implied by the probable capture in the 3:1 resonance with Umbriel. 

In a same time, \citet{peale1988b} shows that Miranda and Ariel are too small for being only heated by tidal effect. He gives the idea to consider another phenomenon like a catastrophic event to increase eccentricities. He also proposes to introduce reliable rheology models and thermal parameters.

%

Here, we study the possible increase of eccentricities by orbital resonance and the resulting endogenic heating by a coupled orbital thermal evolution model. The thermal evolution considers  radiogenic and tidal heating involved by the change 
in the orbital elements when a pair of satellites passes through a mean-motion resonance. The orbital evolution modelizes this passage through the resonance with an averaged 3 body problem, and the coupling of the two parts
depends on the tidal ratio $(k_2/Q)_s$ for the satellite, computed in the thermal module and used in the dynamical one. We consider the heating on each satellite involved. By different rheological models (Maxwell, Burgers and Andrade), 
we compute the rigidity and the viscosity depending on the temperature inside the satellite. We implement next the ratios $(k_2/Q)_s$  for each satellite and 
make evolve the orbital elements with these new ratios. We insist on the fact that we propose a coupled model. We do not consider independently the orbital and thermal evolutions but exchange information
between both modules during the whole simulation.

We present the thermal evolution in the first section containing the process of resolution of the heat equation with source terms for a one dimensional sphere. These source terms are detailed and 
the thermal parameters are defined according to a homogenous mixture of silicates and ices. We also introduce the three rheological models and the associated computation of the dissipation function $Q$.
The second section presents the dynamical module which introduces the resonance 3:1 between Miranda and Umbriel. The Hamiltonian formalism is used to model the 3 body problem and its averaging.
The resulting dynamics is successfully compared with numerical outcomes of the complete 3 body problem. The coupled thermal and orbital evolution is presented in the third section and applied to the pair of satellites Miranda-Umbriel
in the fourth section. This latter is divided in two cases : a {\it nominal scenario} showing the coupled evolution of Miranda with realistic thermal and dynamical parameters/variables, and an {\it extreme scenario} considering higher
orbital eccentricity to try to enhance tidal heating. These two cases show the difficulty to heat a satellite starting form a relatively cold initial state with uniform interior temperatures and surface temperatures in radiative 
thermal equilibrium. Finally, we present in the last section the conclusions and perspectives.

\section{Thermal Evolution}
\label{Thermal_Model}

The satellites of the outer Solar System have various components. Although some of them are exclusively composed by rocks (i.e. Io), the majority of the moons are composed of silicates and ices. This mixture is sometimes homogenous or forms several layers creating a differentiated satellite. Constraints on the composition are given by methods like infrared spectrophotometry which shows for Miranda, a surface composed mainly by water ice \citep{brown1984}. Miranda's bulk density $\rho$=1200  kg/m$^3$  suggests that its interior is mainly composed of water ice and silicate rocks. Whether it is differentiated or not is not known because of the lack of sufficient information on external gravitational field. 

Following the accretion from a mixture of rock and ice, Miranda could have started differentiating if there was sufficient internal heating. The diverse and exotic surface and coronae suggest upwellings of warm material below the surface.  The relatively young age and geology of the coronae is consistent with a temporary geological activity after its formation. It is likely that this internal activity was not active long time enough to alter the whole surface and differentiate the interior \citep{greenberg1991}.  The timing is uncertain  but  such internal activity  could be caused by tidal heating and explained if Miranda was temporarily in a resonant obit with a forced eccentricity. 

We start the calculations considering Miranda as a homogenous mixture of silicate rocks and water ice.  The phase changes of ice are complex and the involvement of components like methane or ammonium, which decrease the melting point of temperature, complicates the study of internal structures evolution \citep{hussmann2009}.  The suspicion of liquid water in some of the moons \citep{hussmann2006} is only validated  by an internal heating of the satellite.  The main heating sources we consider are due to the radiogenic decay of elements in the silicates and the tides due to gravitational effects.

Considering typical ice and silicate densities of $\rho_i=917$\ kg/m$^3$ and $\rho_s=2500$\ kg/m$^3$,  the mass  and volume fraction of silicate rocks  are x$_s$= 0.37$\%$ and f$_s$= 0.45$\%$  respectively. The specific heat   C$_p$ = 900 J kg$^{-1}$ K$^{-1}$  and thermal conductivity  $k$ = 5.2  W m$^{-1}$ K$^{-1}$ of the mixture are calculated based on  mass and volume fractions of the silicates and ices (cf. Table \ref{paramthermiques}):  
  
\begin{equation*}
\begin{split}
 C_p&=x_s\ C_{ps}+(1-x_s)\ C_{p_{i}}\\
 k&=f_s\ k_s+(1-f_s)\ k_{i}\ ,
\end{split}
\end{equation*}
where the $x_s$ and $f_s$ correspond to the mass and volume fractions of silicates respectively. 

Heat transfer occurs principally by  conduction. The transfer of heat by conduction in a one dimension spherically symmetric body is described by the following differential equation (e.g. \citet{schubert2001}): 

\begin{equation}
 \frac{\partial{T(r,t)}}{\partial{t}}=  \alpha \left[  \frac{ \partial^2  T(r,t) } {\partial r^2} +   \frac{2}{r}    \frac {\partial T(r,t)} {\partial r} \right]  + \frac{H}{\rho C_p}\ ,
\label{heat_equation}
 \end{equation}  

\noindent which shows  variation of temperature as a function of satellite radius ($r$) and time ($t$). The parameter $\alpha = \frac{k}{\rho\ C_p}$  is the  thermal diffusivity.  $H$ is the rate of internal heat generation: in our case, we assume radiogenic heating and  tidal dissipation. The heat transfer problem is solved numerically using the  finite differences. The surface temperature $T_{surf}$ is set to the equilibrium temperature $T_{eq}=84$ K and is kept constant along the simulation. To start the calculations, a constant initial temperature profile is assumed with $T(r)=T_{surf}$.  In the center of the satellite ($r=0$), we assumed thermal symmetry i.e,  $\partial{T(0,t)} /  \partial{r}  =  0 $. \\

\begin{table}[h!]
\caption{Physical parameters of the thermal model.}
\vspace{0.2cm}
\centering
\begin{tabular}{l|cccccc}
\hline
& Symbol & Unit & Ice &  Silicate rocks  & Homogenous body \\
\hline
Density  & $\rho$  & kg  m$^{-3}$ & 917& 2500 & 1200 \\
Specific heat & Cp & J kg$^{-1}$ K$^{-1}$ &  888.7 & 920  & 900 \\
Conduction & k & W m$^{-1}$ K$^{-1}$ &  5.4  & 4.2  & 5.2 \\
Rigidity & $\mu$ & Pa  &  4.5 $\times $10$^9$ & 65 $\times$  10$^9$ & 27 $\times$  10$^9$\\
\hline
 \end{tabular}
\label{paramthermiques}
\end{table}

In the thermal evolution calculations, we did not consider changes in porosity and  the average radius is assumed to remain constant over the simulation time (3-6 Myr).  Note that the characteristic time scale of the conduction is proportional to R$^2$  $\alpha^{-1}$ and is  $\approx$ 360 million years.  

Heat is generated inside the silicate part of the satellites through the radioactive decay of  unstable isotopes. The energy emission and the rate of decay depend on the species of radioactive isotope.  More than 98\% of the total radiogenic heat  arises from the decay of the single isotopes of uranium $^{238}$U, $^{235}$U, of thorium $^{232}$Th and order of $1\%$ for potassium $^{40}$K. In the first stages of the evolution, the short-lived radioactive elements $^{26}$Al, $^{60}$Fe and $^{53}$Mn, have a primordial role but are insignificant later. In this study, we consider the radioactive data for the long-lived radioactive elements described in \citet{douce2011} for the radiogenic heating. The short-lived elements will be used for the initial temperature profile (cf. Section \ref{initialconditions}). These elements are gathered in Table \ref{LLRI} from \citep{douce2011}. 

Taking concentrations consistent with the present Earth's mantle \citep{kargel1993},  the  present day radioactive heat production in the mass fraction of silicate rocks of Miranda is  7 $\times$ $10^{-12}$ W/kg  or  $\approx$ $10^{8}$ W.  With the heat capacity of C$_p$ = 900 J kg$^{-1}$ K$^{-1}$, the rate of increase in temperature due to radioactive decay is  only $\approx$ 0.2 K over one million year.   The short-lived radioactive elements on the other hand can provide   2 $\times$ $10^{-7}$ W/kg  or  $\approx$ 5 $\times 10^{12}$ W over the first few million years.

\begin{table}[h!]
\caption{Decay information for the long and short-lived radioactive elements \citep{douce2011}. The parameters $H_{rad_i}$, $\lambda_i$, $C_i^0$ and $C_i^c$ are respectively the rate of radioactive heat production per kg of the initial parent isotope, the decay constant, the half-lived time and the current and initial (4.56 Ga ago) isotopic abundances  of each element.}
 \vspace{0.2cm}
\centering
 \begin{tabular}{l|ccccc}
\hline
&$H_{rad_i}$&$\lambda_i$&$t_{dv_i}$&$C_i^0$&$C_i^a$\\
Isotope&W kg$^{-1}$&$s^{-1}$&yrs&&\\
\hline
$^{238}$U&$9.46\times 10^{-5}$&$4.19\times10^{-18}$&$4.47\times10^{9}$&$0.992\ 75$&\\
$^{235}$U&$5.69\times 10^{-4}$&$3.12\times10^{-17}$&$7.04\times10^{8}$&$0.007\ 20$&\\
$^{232}$Th&$2.64\times10^{-5}$&$1.56\times10^{-18}$&$1.41\times10^{10}$&$1$&\\
$^{40}$K&$2.92\times10^{-5}$&$1.72\times10^{-17}$&$1.28\times10^{9}$&$1.17\times10^{-4}$&\\
\hline
\hline
$^{26}$Al&$4.55\times10^{-1}$&$3.06\times10^{-14}$&$7.17\times10^{5}$&$0$&$5.8\times10^{-5}$\\
$^{60}$Fe&$7.19\times10^{-2}$&$1.46\times10^{-14}$&$1.50\times10^{6}$&$0$&$7\times10^{-7}$\\
$^{53}$Mn&$6.38\times10^{-3}$&$5.87\times10^{-15}$&$3.74\times10^{6}$&$0$&$9\times 10^{-6}$\\
\hline
 \end{tabular}
\label{LLRI}
\end{table}



Tidal dissipation may produce enough heat to keep the internal temperatures, depending on the orbital eccentricity as well as the internal structure and the rheology.  The quantity that characterizes the global dissipation resulting from the non-elastic rheology is the quality factor Q, defined as the ratio of the dissipated energy during one cycle of sinusoidal straining, to the peak energy stored in the system.\\  


For a homogeneous spherical incompressible body with surface gravity g and density $\rho$, the surface potential Love numbers of degree 2, is expressed as:
\begin{equation}
k_2=\frac{3}{2}\ \bigg(1+\frac{19 \tilde \mu}{2\ \rho\ g\ R}\bigg)^{-1}\ ,
\end{equation}

\noindent where $R$ is the radius of the body and $\tilde \mu$ is a complex rigidity obtained by applying the correspondence principle \citep{peltier1974}.  Its expression for different rheological models is given in this section. Among these models Maxwell rheology provides the simplest non-elastic phenomenological rheology adequate for describing dissipation occurring during tidal forcing.  The stress relaxation behavior is described in terms of the Maxwell time $\tau_M =\eta / \mu $. For forcing periods less than the characteristic Maxwell time $\tau_M$,  $ t < \tau_M$ the elastic response predominates and $\tilde \mu \approx \mu $.  The dissipative effects are negligible. For much longer  forcing periods $ t > \tau_M$ the viscous response predominates  and the material behaves like a fluid, $\tilde \mu \approx 0 $.   Maxwell relaxation time for icy moons are in the order of days with a viscosity of  $\eta = 10^{15}$  Pa s,  and rigidity  $\mu = 4.5 \times 10^{9}$  Pa .


There is  little known about the exact rheology parameters  of outer Solar System satellites.  Their elastic properties can be estimated using the Voigt-Reuss-Hill approximation which provides an arithmetic mean between Voigt  and Reuss models \citep{mavko2009}:

\begin{equation}
 \mu_{Voigt}=x_{s}\mu_{s}+(1-x_{s})\mu_{i}\ ,
\end{equation}
and Reuss rigidity :
\begin{equation}
 \mu_{Reuss}= \bigg[\frac{x_{s}}{\mu_{s}} + \frac{(1-x_{s})}{\mu_{i}}\bigg]^{-1}\ ,
\end{equation}
where $\mu_s$ and $\mu_i$ are the rigidities of silicates and ices respectively (cf. Table \ref{paramthermiques}).

Rheology of ice can be complicated, involving several different deformation mechanisms, some of which are non-Newtonian. We  assume a temperature dependent ice rheology, with a Newtonian viscosity $\eta(T)$ that takes the form (see e.g. \citet{parmentier2007}):

\begin{equation}
\eta (T)=\eta_0 \ exp \left[\frac{E_a}{R_gT_m} \bigg(\frac{T_m}{T}-1\bigg)\right]\ ,
\label{Newton_viscosity} 
\end{equation}

where $T_m$ is the  reference temperature and $\eta_0$ the viscosity at  $T=T_m$. The constants $E_a=50\ 10^3$J/mol and $R_g=8.31$ J/mol/K  are the activation energy and the gas constant respectively.

The tidal deformation and resulting deformation of the  ice-rock mixture can be calculated using rheological models which combine elastic \citep{karato1998} and viscous deformations.  We consider in this study three rheological models: Maxwell, Burgers and Andrade. 



The linear viscoelastic Maxwell rheological model gives the complex rigidity $\tilde{\mu}=Re(\tilde{\mu})+Im(\tilde{\mu})$ as:
\begin{equation}
 \tilde{\mu}(\omega)=\frac{\mu\ \eta^2\ \omega^2}{\mu+\eta^2\omega^2}+i\ \frac{\mu^2\ \eta\ \omega}{\mu+\eta^2\omega^2}\ ,
\end{equation}
where the tidal forcing frequency $\omega$ equals the mean motion $n$ of a synchronously rotating satellite; $\mu$ and $\eta$ are the elastic rigidity and the 
steady-state viscosity respectively. 

The Maxwell model tends to overvalue the elastic response of bodies, associated with high viscosities. However, the model depends only on 2 parameters which constitute a big advantage compared to other models. . 


 \begin{figure}[h!]
\begin{center}
\includegraphics[scale=0.4]{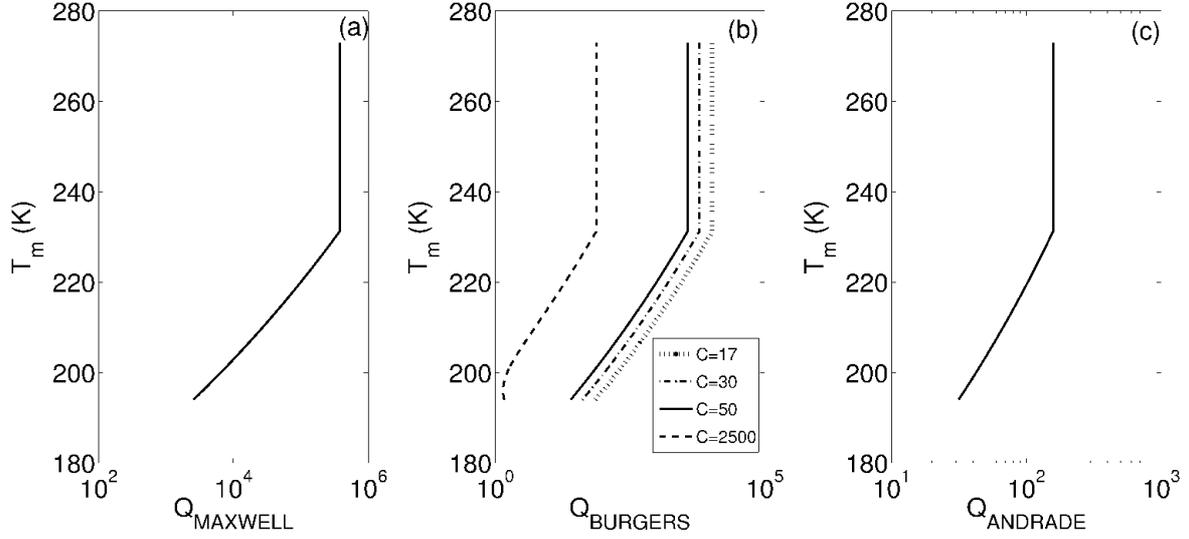}
\caption{Value of $Q$ obtained with Maxwell (a), Burgers (b) and Andrade (c) models in function of the melting temperature $T_m$.}
\end{center}
\label{dissipationQ}
\end{figure}

The Burgers rheology considers a long and a short-term viscosity and is therefore more generic than the Maxwell model \citep{karato1998}:
\begin{equation}
\tilde{\mu}(\omega)=\frac{\omega^2 \left( C_1 -  \eta_{1} C_2 /  \mu_{1} \right)}{ C_2^2+ \omega^2 C_1^2} + i\ \frac{\omega \left( C_2  + \eta_{1}  \omega^2  C_1 /  \mu_{1} \right)}{ C_2^2+ \omega^2 C_1^2}  \ ,
\end{equation}
with
\begin{eqnarray}
C_1 = \frac{1}{\mu_1} +\frac{\eta_1}{\mu_1 \eta_2} +\frac{1}{\mu_2} \\
C_2 = \frac{1}{\eta_2} - \frac{\eta_1}{\mu_1 \mu_2} \omega^2 \ .
\end{eqnarray}
The Burgers model is more efficient than the Maxwell model in the case for instance of the response of terrestrial glaciers to tidal forces \citep{reeh2003}. In icy satellite research, the Burgers body has been applied 
to calculate the despinning of Iapetus \citep{robuchon2010} and the tidal response of Enceladus \citep{shoji2013}. Like \citet{shoji2013}, we assume $\mu_2=\mu_1$ and 
vary $\eta_2/\eta_1$ between 17 and 50.  As upper limit we can also consider  $\eta_2/\eta_1=2500$ as in \citep{shoji2013}. The Burger rheology is relatively more complex than Maxwell model to  incorporate since it requires adjustment of  4 parameters.  

Andrade rheology is an empiric model based on model of viscous fluid in metals \citep{andrade1910}. Resumed by \citep{efroimsky2012} in the case of bodies close to spin-orbit resonances, the model is given by :
\begin{equation}
\tilde{\mu}= \frac{1}{\mu}   + \omega^{-\alpha}  \beta  \cos{\frac{\alpha \pi}{2}} \ \Gamma (\alpha+1) -i \frac{1}{\eta \omega}   - \omega^{-\alpha}  \beta  \sin{\frac{\alpha \pi}{2}} \ \Gamma (\alpha+1) \ ,
\end{equation}
where the parameter $\alpha=0.33\ (0.3-0.38)$ is fixed like for Enceladus \citep{rambaux2010},
\begin{equation}
 \beta = \mu^{\alpha-1} /  \eta^{\alpha}  \approx \bigg[1 \times 10^{-13}; 1 \times 10^{-11}\bigg]\ ,
\end{equation}
and  $\Gamma$ is the gamma function. The number of parameters to be determined in this model makes its handling complicated and difficult compared to Maxwell and Burger models. However, unlike the Maxwell model, the Andrade model can account for the ice anelastic response when forced at periods smaller than the material's Maxwell time (Efroimsky, 2012). The three  rheology models are compared in Figure \ref{dissipationQ}.\\

For a synchronously rotating body in an eccentric orbit  the rate of energy dissipation is (see e.g \citet{peale1999}):
\begin{equation}
\frac{dE}{dt}= \bigg(\frac{k_2}{Q}\bigg)_{\!\!s}\ \frac{\mathcal{G}M^2\ n\ R_s^5}{a^6}\bigg(\frac{21}{2}\ e^2+\frac{3}{2}\ \sin^2 \epsilon\bigg)\ ,
\label{energydissipation}
\end{equation}

\noindent where $\mathcal{G}$ is the gravitational constant, $M$ the mass of the planet and $R_s$ is the satellite radius.
The elements $a$, $e$, $\epsilon$ are the semi-major axis, the eccentricity and the obliquity of the satellite, $n$ is the mean motion. 
Their mean values are given following the JPL website. The value of the radius for the planet $R_p$ corresponds to the values of $J_2$ and $J_4$.  The physical parameters considered for Uranus are in Table \ref{paramphysiqueUranus}.
The physical parameters and the orbital elements of the satellites are gathered in Tables \ref{satellitesphysicalproperties} and \ref{elementsorbitaux} respectively.

 \begin{table}[h!]
\caption{Uranus' physical properties}
\vspace{0.2cm}
 \centering
  \begin{tabular}{l|cc}
   \hline
 Parameter (unit)&Value&Reference\\
 \hline
 GM (km$^3$/s$^2$)&$5\ 793\ 964\ \pm 6$&\citet{jacobson2007}\\
 $J_2\ \times 10^6$&$3\ 341.29\pm0.72$&\citet{jacobson2007}\\
$J_4\ \times 10^6$&$-30.44\pm 1.02$&\citet{jacobson2007}\\
 \hline
  \end{tabular}
\label{paramphysiqueUranus}
 \end{table}

The derivation of the formula (\ref{energydissipation}) assumes that the body is incompressible, the rotation is uniform and synchronous.

The dissipated energy is associated with orbital parameters as it arises from two distinct sources of time dependence in the tide: time variation in the distance to the tide-raising planet, and the optical libration (the relative rocking motion of a uniformly rotating satellite relative to the planet that results from the nonuniform motion in the elliptic orbit).  The dissipation depends on semi-major axis, eccentricity and inclination of the orbit, these last two parameters varying with the encounter of resonances. Although the effect is not significant, we keep the effect on the obliquity which is computed at the Cassini State 1 by \citet{noyelles2010}: 

\begin{equation}
 \epsilon_{eq}\approx \frac{\sin I}{\alpha/\dot{\Omega}+\cos I}\ ,
\end{equation}
where
\begin{equation}
 \alpha=\frac{3}{2}\ \frac{(C-A)n}{C}\ ,
\end{equation}
where $A$ and $C$, the principal moments of inertia, are given by:
\begin{equation}
 \begin{split}
  A&=\frac{4}{15}\ \rho\ \pi\ abc\ (b^2+c^2)\\
C&=\frac{4}{15}\ \rho\ \pi\ abc\ (a^2+b^2)\ ,
 \end{split}
\end{equation}
where $a$, $b$ and $c$ depend on the satellite shape. Their values are resumed in Table \ref{satellitesphysicalproperties}.

 \begin{table}[h!]
\caption{Physical properties of Miranda and Umbriel}
\vspace{0.2cm}
 \centering
  \begin{tabular}{l|ccc}
   \hline
 Parameter (unit)&Miranda&Umbriel&Reference\\
 \hline
 GM (km$^3$/s$^2$)&$4.4 \pm 0.4$&$81.5\pm5.0$&\citet{jacobson2007}\\
 Mean Radius (km)&$235.8\pm0.7$&$584.7\pm2.8$&\citet{thomas1988}\\
 Subplanetary equatorial &&&\\
 radius (km)&$240.4$&&\citet{Archinal2009}\\
 Along orbit equatorial &&&\\
 radius (km)&$234.2$&&\citet{Archinal2009}\\
 Polar radius (km)&$232.9$&&\citet{Archinal2009}\\
 \hline
  \end{tabular}
 \label{satellitesphysicalproperties}
 \end{table}

 \begin{table}[h!]
 \caption{Mean orbital elements of the five main satellites at J2000 \citep{laskar1987}: $a$ is the semimajor axis, $e$ the eccentricity, $\omega$ the pericenter, $M$ the mean anomaly, $I$ the inclination, $\Omega$ the ascending node, 
   $n$ the mean motion. The variables $P$ and $P_{\Omega}$ stand for the orbital and the node periods respectively.}
  \vspace{0.2cm}
 \centering
  \begin{tabular}{l|ccccccccc}
  \hline
   Satellites & $a$ & $e$ & $\omega$ & $M$ & $I$ & $\Omega$ & $n$ & $P$& $P_{\Omega}$\\
      & (km)&& (deg)& (deg)& (deg)& (deg)&(deg/day)&(days)&(yr)\\
  \hline
  Miranda  & $129\ 900$ & $0.0013$ &  $68.312$ & $311.330$ & $4.338$ & $326.438$ & $254.6906576$ & $2.520$ & $17.727$ \\
  Umbriel  & $266\ 000$ & $0.0039$ &  $84.709$ &  $12.469$ & $0.128$ &  $33.485$ &  $86.8688879$ & $8.706$ & $126.951$ \\
  \hline
  \end{tabular}
 \label{elementsorbitaux}
 \end{table}


Internal heating due to tidal dissipation would increase the internal temperatures. The temperature dependent viscosity decreases with increasing temperatures resulting in an increase of $k_2/Q$ and tidal dissipation. The orbital and thermal evolution are coupled through the parameter $k_2/Q$  which affects the orbital parameters and resonances as described in the section \ref{Results}.

\section{Dynamical Model}
\label{Dynamical_Model}


In this section, we introduce the modelization of the dynamical problem. The N-body problem of Uranus and its five main satellites has already been studied in details in \citep{verheylewegen2013}
where we consider a planetocentric reference frame and
\begin{itemize}
 \item the gravitational interactions of the five main satellites seen as point masses, 
 \item the oblateness of Uranus up to the second order $J_2$ and $J_4$ (cf. Table \ref{paramphysiqueUranus}).
\end{itemize}
To make evolve the system in time, we also add to these latest perturbations, the tidal effect using the Kaula's formulations (see e.g. \citet{Yoder1981}):
\begin{equation}
\begin{split}
 \frac{da}{dt}&=3\ \bigg(\frac{k_2}{Q}\bigg)_{\! p}\ \frac{n\ m\ R_p^5}{a^4\ M}\ \bigg(1+\frac{51}{4}e^2\bigg)-21\ \bigg(\frac{k_2}{Q}\bigg)_{\! s}\ \frac{n\ MR_s^5}{a^4m}\ e^2 \\
 \frac{de}{dt}&=\frac{57}{8}\ \bigg(\frac{k_2}{Q}\bigg)_{\! p} \frac{n\ m}{M}\bigg(\frac{R_p}{a}\bigg)^5\ e-\frac{21}{2}\ \bigg(\frac{k_2}{Q}\bigg)_{\! s}\ \frac{n\ M}{m}\bigg(\frac{R_s}{a}\bigg)^5e\ ,
\end{split}
\label{Kaula}
\end{equation}
where the index $p$ and $s$ refer to the planet and the satellite respectively, $R_p$ is the mean radius of the planet and $m$ the mass of the satellite (cf. Tables \ref{paramphysiqueUranus} and \ref{satellitesphysicalproperties}). 
Due to the small oblateness of Uranus, the resonances overlap and the assumption of an isolated resonance holds only in the particular cases 
of small inclinations and eccentricities. Therefore we have to take into account the six resonant arguments of second order in the 3:1 mean-motion resonance between Miranda and Umbriel, which are:  
\begin{center}
 \begin{tabular}{@{}ll@{}}
 $2\ \theta_1=\lambda_5-3\lambda_2+2\Omega_5$&\quad\quad $[I_M^2]$   \\
 $2\ \theta_2=\lambda_5-3\lambda_2+\Omega_5+\Omega_2$&\quad\quad $[I_MI_U]$  \\
 $2\ \theta_3=\lambda_5-3\lambda_2+2\Omega_2$&\quad\quad $[I_U^2]$ \\
 $2\ \theta_4=\lambda_5-3\lambda_2+2\varpi_2$&\quad\quad $[e_U^2]$ \\
 $2\ \theta_5=\lambda_5-3\lambda_2+\varpi_5+\varpi_2$&\quad\quad $[e_Me_U]$ \\
 $2\ \theta_6=\lambda_5-3\lambda_2+2\varpi_5$&\quad\quad $[e_M^2]$  \ .
\end{tabular}
\label{Resonant_Argument}
\end{center}
\vspace{0.4cm}
\noindent
where, in the left column, $\theta_k,\ k=1:6$ are the resonant arguments for the primary resonances with $\lambda_i$, the mean longitudes, $\Omega_i$ the ascending nodes and, $\varpi_i$ the longitudes of the pericenters. 
The index $5$ and $2$ stand respectively for Miranda and Umbriel, following the label given by chronological order of discovery of each satellite. The right column is the type of the resonance and corresponds to the first non-zero term associated with the cosine of the angle $\theta_i$ 
in the perturbative potential.

With new powerful methods, we studied the N-body problem and an averaged form in \citet{verheylewegen2013}. We retrieved the main result of \citet{tittemore1989,tittemore1990, malhotradermott1990}, namely the high value of Miranda due to the capture in the 3:1 mean-motion resonance with Umbriel. 
In particular, we showed with the frequency analysis tool \citep{laskar1993} that the exit at $4.5^\circ$ for the inclination of Miranda can be due to the disruption of the primary resonance by a 2:1 secondary resonance, subsequent the capture 
into a 3:1 secondary resonance between the frequency of the libration argument of $\theta_1$ and the frequency of the circulation argument of $\theta_2$ \citep{malhotradermott1990, verheylewegen2013}.\

\par For the dynamical part, our first coupled model was based on the full N-body problem with an artificial increase of the tidal ratio $(k_2/Q)_p$ of the planet: strengthening this ratio allows the increase of the rate of evolution of the system 
providing that the variations of the orbital elements remain adiabatic (see e.g. \citep{malhotra1991}).

\par Since the characteristic timescale for the evolution of the temperature inside the satellites is long and to preserve the physical significations related to the evolution of temperature, 
we choose to develop an averaged form of the full dynamical problem consistent with the idea of a coupling with a thermal long-term evolution and preserve a ratio $(k_2/Q)_p$ in accordance with
the studies of \citet{tittemore1988,tittemore1989,tittemore1990}.

\par Averaged models of the Uranian system have already been studied by different authors 20 years ago.
\citet{tittemore1988} constructed an Hamiltonian in a planar eccentric case and extended it in a inclined circular case in \citep{tittemore1989} considering the 2 body gravitational interaction, the perturbation due to the oblateness
of the planet, the resonant terms, and finally the perturbation due to the secular interactions between the satellites. They eventually obtained an Hamiltonian in canonical coordinates with four degrees of freedom by the addition 
of the two previous cases. \citet{malhotradermott1990} worked with an inclined circular or a planar eccentric Hamiltonian separately to analyze the role of the secondary resonances in the 3:1 mean-motion inclination or eccentricity resonances respectively.

\par We choose to implement the method explained in the case of the Saturnian System in \citet{champenois1998}, to select rigorously the terms needed in our modelisation and to obtain an averaged Hamiltonian depending on the inclinations 
and on the eccentricities in the same time. The reason of this choice is the following : the thermal heating is more efficient when we have an increase in eccentricities (cf. Equation \ref{energydissipation}) but we also think that the increase in inclination 
for Miranda is a key point of the evolution of the system and that the capture into the resonance $\theta_1$ is necessary to have a good approach of the problem.

\subsection{The Averaged Hamiltonian}

By introducing Jacobian coordinates, the usual Hamiltonian is written to the first order on satellite masses \citep{tittemore1988}:
\begin{eqnarray}
 \mathcal{H}&=&-\sum_{i=1}^{N} \frac{\mathcal{G}Mm_i}{2 a_i}\left[1+\sum_{n=1}^2 J_{2n} \bigg(\frac{R_p}{a_i}\bigg)^{2n} P_{2n} (\sin{\phi_i})\right]\nonumber \\
&-&\mathcal{R}\ ,  
\label{Hamilt_init}
\end{eqnarray}
\noindent
where $m_i$ is the mass of the satellite $i$. The variable $R_p$ is the radius of the planet corresponding to the values of $J_2$ and $J_4$, $a_i$ is the semi-major axis of the satellite $i$. 
The $N$ first terms in the Hamiltonian~(\ref{Hamilt_init}) consider
the two-body interaction between Uranus and each satellite. The second ones are the perturbation due to the oblateness of the planet developed in classical
Legendre polynomial, with $\phi_i$ the latitude of the satellite $i$. For the Uranian system we only consider the known spherical harmonics $J_2$ and $J_4$.
The last term is the perturbation due to the third body contained in the disturbing function $\mathcal{R}$, written in the first order of masses as (e.g. \citet{champenois1998}):
\begin{equation}
 \mathcal{R}_{ij}=\frac{1}{a_j}\mathcal{G}Mm_j\ \bigg(\frac{a_j}{r_{ij}}-a_j\frac{\mathbf{r_i}\cdot\mathbf{r_j}}{r_j^3}\bigg)
\label{disturbingfunction_ext}
\end{equation}
 for the outer perturbation by a satellite $j$ on a satellite $i$ and 
\begin{equation}
 \mathcal{R}_{ji}=\frac{1}{a_j}\mathcal{G}Mm_i\ \bigg(\frac{a_j}{r_{ij}}-a_j\frac{\mathbf{r_i}\cdot\mathbf{r_j}}{r_j^3}\bigg)
\label{disturbingfunction_int}
\end{equation}
in the case of the inner perturbation by a satellite $i$ on a satellite $j$. \\ 

\par Considering the effect of Umbriel on Miranda (resp. of Miranda on Umbriel), we can rewrite the external disturbing function (\ref{disturbingfunction_ext})
(resp. the internal disturbing function (\ref{disturbingfunction_int})) as:

 \begin{equation}
  \mathcal{R}_{52}=\frac{1}{a_2}\ \mathcal{G}Mm_2\ \bigg(\frac{a_2}{\Delta_{52}}-a_2\ \frac{\mathbf{r_5}\cdot\mathbf{r_2}}{r_2^3}\bigg)
\label{perturbativefunction_initMU}
\end{equation}
\begin{equation}
  \mathcal{R}_{25}=\frac{1}{a_2}\ \mathcal{G}Mm_5\  \bigg(\frac{a_2}{\Delta_{52}}-a_2\ \frac{\mathbf{r_5}\cdot\mathbf{r_2}}{r_2^3}\bigg)
\label{perturbativefunction_initUM}
 \end{equation}
where $\Delta_{52}$ is the distance between Miranda and Umbriel. Classically (see e.g. \citet{murray1999}), we expand these perturbative functions
to the second order in eccentricity-inclination and select the long period terms. These period terms are typically of about 100 years. The selection of these
terms are introduced in the following subsections with the objective to determine the perturbative function needed in the averaged model.

\subsubsection{The resonant terms}

The resonant terms are the arguments associated with the second order mean-motion resonance 3:1 between Miranda and Umbriel. Typically we consider
the six possible resonant angles summarized in the beginning of the section~\ref{Dynamical_Model}. Each resonant angle in the perturbative function 
is associated with a Laplace coefficient function $f_k(\alpha),\ k=1:6$ (see e.g. \citet{murray1999}):

\begin{center}
 \begin{tabular}{@{}lcl@{}}
$f_1(\alpha)$&$=$&$\frac{1}{2}\ \gamma_5^2\ \alpha\ b_{3/2}^{(2)}$\\
$f_2(\alpha)$&$=$&$-\gamma_5\ \gamma_2\ \alpha\ b_{3/2}^{(2)}$\\
$f_3(\alpha)$&$=$&$\frac{1}{2}\gamma_2^2\ \alpha\ b_{3/2}^{(2)}$\\
$f_4(\alpha)$&$=$&$\frac{1}{8}e_2^2\ \bigg(17+10\ \alpha\ D+\alpha^2\ D^2\bigg)\ b_{1/2}^{(1)}$\\
$f_5(\alpha)$&$=$&$-\frac{1}{4}e_5\ e_2\ \bigg(20+10\ \alpha\ D+\alpha^2\ D^2\bigg)\ b_{1/2}^{(2)}$\\
$f_6(\alpha)$&$=$&$\frac{1}{8}e_5^2\ \bigg(21+10\ \alpha\ D+\alpha^2\ D^2\bigg)\ b_{1/2}^{(3)}$\ ,
\end{tabular}
\end{center}
where $\gamma_i=\sin \frac{I_i}{2}$, $\alpha=a_5/a_2$ is the ratio of semi-major axes and $b_{s}^{(j)} (\alpha)$ are the coefficients of Laplace defined by (see e.g. \citet{murray1999}):
\begin{equation}
 \frac{1}{2}\ b_s^j(\alpha)=\frac{1}{2\pi}\ \int_0^{2\pi}\ \frac{\cos j\ \Psi\ d\Psi}{(1-2\alpha\cos{\Psi}+\alpha^2)^s}\ ,
\end{equation}

and $D$, $D^2$ are the differential operators related to $\alpha$ of first and second order of these coefficients.\\
 
The six resonant terms selected in this section are part of the {\bf direct} resonant perturbative function. Due to our choice of a planetocentric frame, we need 
to consider also the {\bf indirect} resonant terms (which are absent in the case of a barycentric frame) expressed as:
\begin{eqnarray}
 \mathcal{R}_E&=-\frac{27}{8}\ e_2^2\cos(\lambda_5-3\lambda_2+2\varpi_2)\\
\mathcal{R}_I&=-\frac{3}{8}\ e_2^2\cos(\lambda_5-3\lambda_2+2\varpi_2)\ ,
\end{eqnarray}
for an outer or an inner perturbation respectively.\\

The expressions of the perturbative functions (\ref{perturbativefunction_initMU}) and (\ref{perturbativefunction_initUM})
for the resonant terms are given by:
\begin{eqnarray}
 \mathcal{R}_{52}^R&=&\frac{\mathcal{G}M m_2}{a_2} \ \bigg(\sum_{k=1}^6 f_k\cos{2\theta_k}-a_2\ \mathcal{R}_E\bigg)\\
 \mathcal{R}_{25}^R&=&\frac{\mathcal{G}M m_5}{a_2} \ \bigg(\sum_{k=1}^6 f_k\cos{2\theta_k}-a_2\ \mathcal{R}_I\bigg)\ .
\end{eqnarray}

\subsubsection{The secular terms}
The secular terms are the independent expressions and the terms depending on the difference between nodes and pericenters in 
the expansion of the perturbative functions (\ref{perturbativefunction_initMU}) and (\ref{perturbativefunction_initUM}), each term being
associated with a Laplace coefficient function $C_k(\alpha),\ k=0:4$ (see e.g. \citet{murray1999}):
\begin{center}
 \begin{tabular}{@{}lcl@{}}
$C_0(\alpha)$&=&$\frac{1}{2}\ b_{1/2}^{(0)}$\\
$C_1(\alpha)$&$=$&$\frac{1}{8}\ (2\alpha\ D+\alpha^2 D^2)\ b_{1/2}^{(0)}$\\
$C_2(\alpha)$&$=$&$-\frac{1}{2}\ \alpha\ b_{3/2}^{(1)}$\\
$C_3(\alpha)$&$=$&$\frac{1}{4}\ (2-2\alpha D-\alpha^2\ D^2)\ b_{1/2}^{(1)}$\\
$C_4(\alpha)$&$=$&$\alpha\ b_{3/2}^{(1)}$\ .
\end{tabular}
\end{center}
The expressions of the perturbative functions (\ref{perturbativefunction_initMU}) and (\ref{perturbativefunction_initUM}) for the secular terms are given by:
\begin{eqnarray}
 \mathcal{R}_{52}^S&=&\frac{\mathcal{G}M m_2}{a_2} \ \bigg(C_0+C_1\ (e_5^2+e_2^2)+C_2\ (\gamma_5^2+\gamma_2^2)+C_3\ e_5\ e_2\cos(\varpi_2-\varpi_5)+C_4\ \gamma_5\ \gamma_2\cos(\Omega_2-\Omega_5)\bigg)\\
 \mathcal{R}_{25}^S&=&\frac{\mathcal{G}M m_5}{a_2} \ \bigg(C_0+C_1\ (e_5^2+e_2^2)+C_2\ (\gamma_5^2+\gamma_2^2)+C_3\ e_5\ e_2\cos(\varpi_2-\varpi_5)+C_4\ \gamma_5\ \gamma_2\cos(\Omega_2-\Omega_5)\bigg)\ .
\end{eqnarray}
\subsubsection{The oblateness terms}

It remains to complete the perturbative function by the oblateness term. An averaging version to the second order of this term can be find in \citet{murray1999}:
\begin{equation}
\mathcal{R}_{i}^A= \frac{\mathcal{G}M}{2a_i}\ \bigg[\frac{3}{2}\ J_2\ \bigg(\frac{R_p}{a_i}\bigg)^2-\frac{9}{8}J_2^2\ \bigg(\frac{R_p}{a_i}\bigg)^4-\frac{15}{4}
J_4\ \bigg(\frac{R_p}{a_i}\bigg)^4\bigg]\ e_i^2\\
-\frac{\mathcal{G}M}{2a_i}\ \bigg[\frac{3}{2}\ J_2\ \bigg(\frac{R_p}{a_i}\bigg)^2-\frac{27}{8}J_2^2\ \bigg(\frac{R_p}{a_i}\bigg)^4-\frac{15}{4}
J_4\ \bigg(\frac{R_p}{a_i}\bigg)^4\bigg]\ \sin^2I_i\ .
\end{equation}

\subsection{The equations of motion}
The equations defined up to now depend on the orbital elements of the satellites $(a,e,i,\varpi,\Omega)$ and on the resonant arguments $\theta_k, \ k=1:6$.
We choose to work here with the Lagrangian variables (e.g. \citet{duriez1977}):
\begin{eqnarray}
 z_i    &=&e_i\ \exp\ (\sqrt{-1}\ \varpi_i)\\
 \zeta_i&=&\gamma_i\ \exp(\sqrt{-1}\ \Omega_i)\ ,
\label{eq_Lagrange}
\end{eqnarray}
with $i=2,5$ for Umbriel or Miranda. The definition of these variables avoids the problem of indetermination of the pericenters and/or the nodes when the eccentricities and/or the inclinations are equal to zero. In these 
variables and following \citet{duriez1977}, the equations of the motion are written as:

\begin{eqnarray}
\frac{da_i}{dt}&=&\frac{2}{n_ia_i}\ \frac{\partial \mathcal{R}_i^L}{\partial \lambda_i}\\
\frac{dz_i}{dt}&=&\frac{\sqrt{-1}\ \phi_i}{n_ia_i^2}\ \bigg[2\frac{\partial \mathcal{R}_i^L}{\partial\bar{z_i}}+\frac{\sqrt{-1}}{(1+\phi_i)}\ z_i\ \frac{\partial \mathcal{R}_i^L}{\partial \lambda_i}+\frac{z_i}{2\phi_i^2}\bigg(\zeta_i\frac{\partial \mathcal{R}_i^L}{\partial \zeta_i}+\bar{\zeta_i}\frac{\partial \mathcal{R}_i^L}{\partial\bar{\zeta_i}}\bigg)\bigg]\\
\frac{d\zeta_i}{dt}&=&\frac{\sqrt{-1}}{2n_ia_i^2\phi_i}\ \bigg[\frac{\partial \mathcal{R}_i^L}{\partial\bar{\zeta_i}}+\sqrt{-1}\ \zeta_i\ \frac{\partial \mathcal{R}_i^L}{\partial \lambda_i}-\zeta_i\ \bigg(z_i\ \frac{\partial \mathcal{R}_i^L}{\partial z_i}-\bar{z_i}\ \frac{\partial \mathcal{R}_i^L}{\partial\bar{z_i}}\bigg)\bigg]\\
\frac{d\lambda_i}{dt}&=&n_i-\frac{2}{n_ia_i}\ \frac{\partial \mathcal{R}_i^L}{\partial a_i}
+\frac{\phi_i}{n_ia_i^2(1+\phi_i)}\ \bigg(z_i\frac{\partial \mathcal{R}_i^L}{\partial z_i}+\bar{z_i}\frac{\partial \mathcal{R}_i^L}{\partial\bar{z_i}}\bigg)+\frac{1}{2n_ia_i^2\phi_i}\bigg(\zeta_i\frac{\partial \mathcal{R}_i^L}{\partial \zeta_i}+\bar{\zeta_i}\ \frac{\partial \mathcal{R}_i^L}{\partial\bar{\zeta_i}}\bigg)\ ,
\end{eqnarray}
with $\phi_i=\sqrt{1-z_i\bar{z_i}}$. The expression $\mathcal{R}_i^L\ i=2,5$ is the perturbative function containing all the long period terms selected in the previous section and is expressed by:
\begin{equation}
 \mathcal{R}_i^L=\mathcal{R}_{ij}^R+\mathcal{R}_{ij}^S+\mathcal{R}_{i}^A\ ,
\end{equation}
with $i=2,5$ for the two satellites Umbriel and Miranda and $j=5,2$ for the perturbation due to the third body. 
We integrate a set of $11$ differential equations: 
\begin{equation}
 \bigg(\frac{da_i}{dt},\frac{dk_i}{dt},\frac{dh_i}{dt},\frac{dq_i}{dt},\frac{dp_i}{dt},\frac{d\Psi}{dt}\bigg)\ ,
\label{set_equations}
\end{equation}
the variables $k_i$, $h_i$, $p_i$ and $q_i$ being defined by
\begin{eqnarray}
 k_i&=&e_i\cos\varpi_i=\mathop{\mathrm{Re}}\ (z_i)\\
 h_i&=&e_i\sin\varpi_i=\mathop{\mathrm{Im}}\ (z_i)\\
 q_i&=&\gamma_i\cos\Omega_i=\mathop{\mathrm{Re}}\ (\zeta_i)\\
 p_i&=&\gamma_i\sin\Omega_i=\mathop{\mathrm{Im}}\ (\zeta_i)\ .
\end{eqnarray}
The variable $\Psi=3\lambda_2-\lambda_5$ is the exact resonant angle.\\

This set of equations of motion (\ref{set_equations}) is integrated with the same Adams-Bashforth-Moulton $10^{\text{th}}$ order predictor-corrector integrator
than in \citet{verheylewegen2013} and we invite the reader to refer to this last article for the validation of the numerical code. To validate the results
of the integration of the averaged model presented here, we proceed in the same way as in \citet{verheylewegen2013}.

\begin{center}
 \begin{figure}[h!]
\includegraphics[scale=0.45]{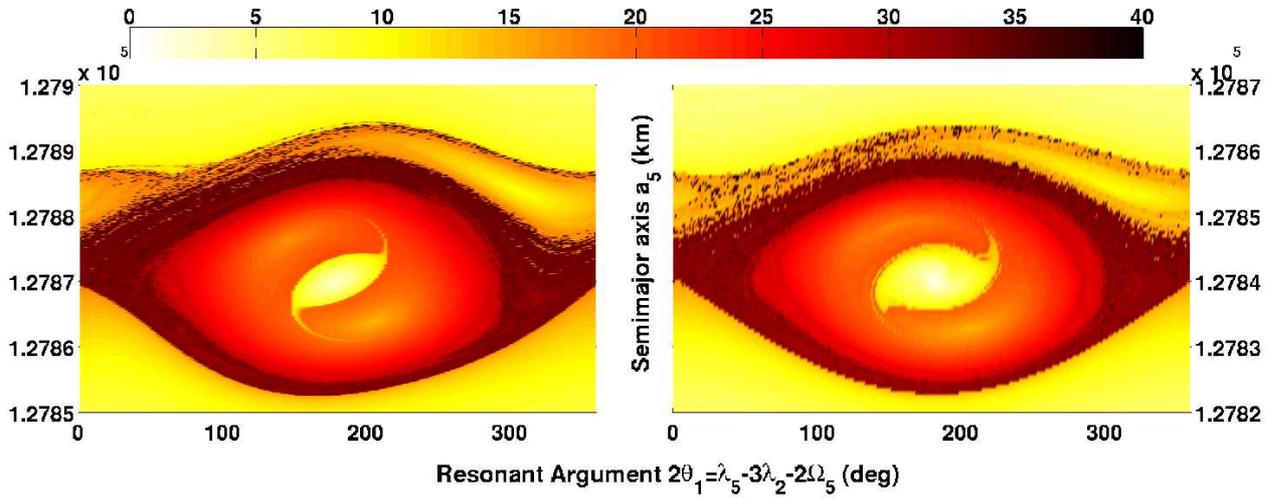}
\caption{Phase spaces semimajor axis $a_5$ versus resonant argument $\theta_1$ resulting from the 3 body problem Uranus, Miranda, Umbriel with the Adams-Bashforth-Moulton integrator over 1500 years (left hand side figure) and
its averaged version (right-hand side figure). For the complete version, the integration step is set to $1/80$ day. The initial conditions are the current ones at J2000 (cf. Tables \ref{paramphysiqueUranus}, 
\ref{satellitesphysicalproperties} and \ref{elementsorbitaux}) except for the mean anomaly $M_5$, the semimajor axis $a_5$
and the inclination of Miranda $I_5$. The two first variables are set respectively between $[0°-360°[$ and $[127850\ {\text km}-127900\ {\text km}]$.
For the averaged version, the integration step is set to $17/300$ years. The semimajor axis $a_5$ is between $[127820\ {\text km}-127870\ {\text km}]$. The initial inclination for Miranda is $4.338^{\circ}$ in both cases.
The colorbar considers the variations of semimajor axis $a_5$ (km) on each simulation.}
\label{eyeofresonance}
\end{figure}
\end{center}

We represent the 3:1 mean motion resonance between Miranda and Umbriel with color maps. To obtain Figure~\ref{eyeofresonance}, we perform $10^4$
numerical integrations over $1500$ years, each one associated with a different initial condition. We choose a time step of $17/300$ years corresponding to $1/300^{\text{th}}$ 
of the smallest nodal period i.e. the nodal period of Miranda. The color scale here contains the variations
of the semi-major axis of Miranda: at the end of the simulation, the color is the difference between the largest value of the semi-major axis and 
the smallest one. In a libration zone, the semi-major axis is locked and we have a difference between the maximum and the minimum value 
near zero, corresponding to a light color. For a more chaotic zone (like a separatrix), the variations are larger, corresponding to a dark color.
This type of color scale has already been compared with the chaos indicator Mean Exponential Growth of Nearby Orbits (MEGNO) \citep{cincotta2000} in \citet{verheylewegen2013}. 

On the left-hand side in Figure \ref{eyeofresonance}, we plot the result of the 3 body problem integration presented in \citet{verheylewegen2013} which
considers the gravitational perturbations between Uranus, Miranda and Umbriel and the oblateness effect and containing in particular all the short period terms. 
On the right-hand side of Figure \ref{eyeofresonance}, we plot the result of the numerical integration of the equations of motion (\ref{set_equations}) containing
the same perturbations as in the complete system but considering only the long period terms. We observe that the two figures are almost the same
with more precision in the case of the integration of the complete system. But our averaged form maintains the 
global points of the dynamics of the mean motion resonance 3:1 between Miranda and Umbriel with the presence of the large separatrix delimiting the 
border of the resonance. We can also distinguish the two secondary resonances zones in the center of libration playing a role in the exit of the 
primary resonance \citep{malhotradermott1990}. We conclude this section by saying that our averaged form is validated by the study of the dynamical part 
of the problem.

\section{Coupling of the dynamical and the thermal parts: procedure}
Coupling studies mixing dynamic and thermal aspects of evolution are rather rare. Indeed they generate difficulties in the combination of the two approaches.
The characteristic times for instance are very different with changes in the orbital motion on several days/years while the temperature inside the satellites
needs millions of years to vary. Therefore, it is sometimes needed to make a choice between a plausible dynamical evolution and a physical thermal evolution.

In \citet{schubert2010}, the authors analyze the role of the resonances in the internal evolution of the satellites: they show the influence of eccentricity pumping on the tidal heating. In particular, for the Galilean satellites, the tidal heating of Io exceeds easily the radiogenic heating involving the well-known volcanism in its surface. An important dynamical element is the Laplace resonance between the three satellites Io, Europa and Ganymede whose libration argument is given by:
\begin{equation}
 \theta=\lambda_1-3\lambda_2+2\lambda_3\ ,
\end{equation}
where $\lambda_i,\ i=1:3$ are the mean longitudes of Io, Europa and Ganymede respectively. This configuration is stable and has the particularity to make evolve Io in a hot state to a cold state and vice-versa when the satellites evolve inside the resonance. This heating cycles process is proposed by \citet{ojakangas1986}. Two important things here are the maintenance of the resonance in time and the convection form of heat transfer for Io. On the same system, we can also cite the work of \citet{showman1997} which presents a coupled dynamical-thermal model for Ganymede in the case of a satellite with homogenous temperature, without radial  variation of temperature. In the case of the Uranian System, there is to our knowledge, no coupled approach. The authors who studied thermal questions in the cases of Miranda and Ariel  (\citealt{dermott1988}; \citealt{peale1988b}; \citealt{tittemore1988,tittemore1989,tittemore1990}) investigated the two aspects separately.\\

The solution procedure of the coupled model is presented in Figure \ref{schema}. In a dynamical point of view, the averaged equations of motion (\ref{set_equations}) are solved and give the orbital parameters. These latter involved in the computation of the tidal heating in the resolution of the heat equation (\ref{heat_equation}). The thermal module provides temperature dependent  viscosity and the value of the dissipation function using Maxwell, Burgers or Andrade models. The new ratio $k_2/Q$ obtained for the satellite is then used in the dynamical module.\\
\begin{center}
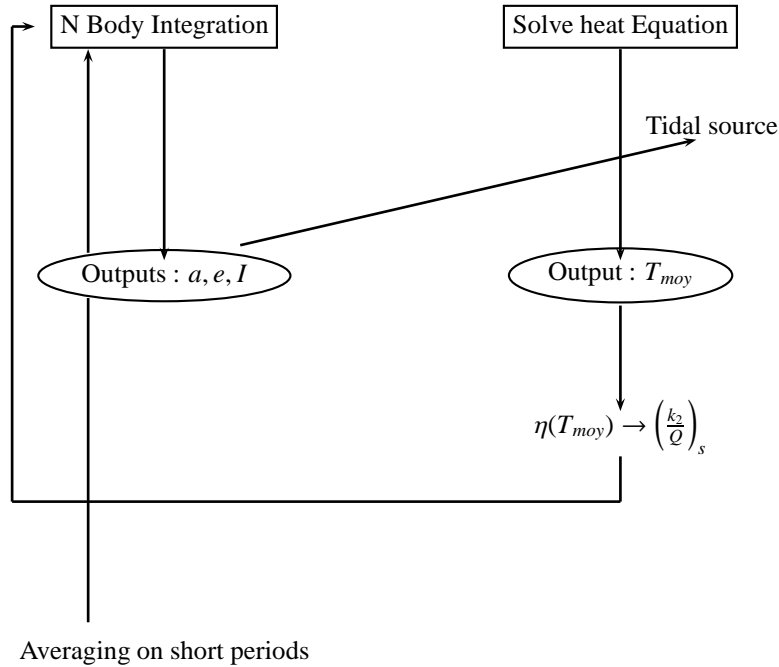
\begin{figure}[h!]
 \begin{center}
\vspace{1cm}
  \begin{pspicture}(0,2)(10,10)
\rput(2,10.3){\psframebox{N Body Integration}}
\rput(8,10.3){\psframebox{Solve heat Equation}}
\psline[linewidth=1pt]{->}(2,10)(2,7.2)
\psline[linewidth=1pt]{->}(8,10)(8,7.2)
\rput(2,7){\psovalbox{Outputs : $a,e,I$}}
\rput(8,7){\psovalbox{Output : $T_{moy}$}}
\psline[linewidth=1pt]{->}(8,6.6)(8,5.2)
\rput(8,5){$\eta(T_{moy})\rightarrow\bigg(\frac{k_2}{Q}\bigg)_s$}
\rput(9.2,9){Tidal source}
\psline[linewidth=1pt]{->}(3,7.4)(9,8.8)
\psline[linewidth=1pt]{-}(8,4.6)(8,4)
\psline[linewidth=1pt]{-}(8,4)(0,4)
\psline[linewidth=1pt]{-}(0,4)(0,10.3)
\psline[linewidth=1pt]{->}(0,10.3)(0.3,10.3)
\rput(2,2){Averaging on short periods}
\psline[linewidth=1pt]{-}(1,2.4)(1,6.7)
\psline[linewidth=1pt]{->}(1,7.3)(1,10)
 \end{pspicture}
\end{center}
\caption{Schematic representation of the coupled approach.}
\label{schema}
\end{figure}
\end{center}

Following its accretion, the satellite's interior starts to cool down since the internal heating dominated by the radiogenic elements decay is not sufficient. If the satellite is captured in a resonance which pumps the  orbital eccentricity sufficiently, the tidal heating becomes important and dominates the radiogenic heating. The viscosity is then decreasing with increasing temperature. If the temperature is high enough, the satellite can be differentiated in several layers with the heavy elements in the core. The tidal dissipation damps orbital eccentricity.  The tidal dissipation in return,  diminishes with lower eccentricities  (cf. Equation \ref{energydissipation}).

\section{Coupling of the dynamical and the thermal parts: results}
\label{Results}

In this section, we will show the difficulty to heat a satellite like Miranda when we start calculations with a relatively cold interior. The coupled simulation focus on a capture in a resonance in eccentricity.
Since there is no geological evidence for resurfacing of Umbriel \citep{smith1986}, we do not consider a rise in eccentricity for Umbriel and we select a resonance of type $e_M^2$ allowing us to observe a larger increase in the eccentricity of Miranda. We suppose that the inclination of Miranda is already $4.5^\circ$. We present the results in different cases with various rheological models.

\subsection{Initial conditions}
\label{initialconditions}

The satellite radius is discretized in $N_i$ points. The temperature is then computed in each predetermined layer $i$.
When the temperature has been calculated in the whole satellite, we compute an averaged temperature 
associated with an averaged viscosity and $(k_2/Q)_s$ for each satellite. The heating by the radiogenic elements is weak because it does not consider the short-lived elements. The order of this heating is $10^{-14}$ W/kg.

As initial conditions, we first consider a homogenous satellite at the constant equilibrium temperature $T_{eq}$, which value is $84$ K for Miranda. 
For comparison, \citet{castillo2007} obtain $T_{eq}=90$ K for Iapetus. To be more realistic we choose another approach consisting of a thermal profile on the entire satellite.
These thermal profiles are determined by the resolution of the equation (\ref{heat_equation}) with the boundary and initial conditions previously defined. As source term, we consider the effect of the radiogenic heating 
since the formation of the satellites over $4.6$\ Gyr. In this case, we have to take into account the short-lived radioactive elements which are listed in the Table \ref{LLRI} \citep{douce2011} because they are active in the early stage 
of the evolution of the satellites.

\begin{figure}[h!]
\begin{center}
\includegraphics[width=120mm]{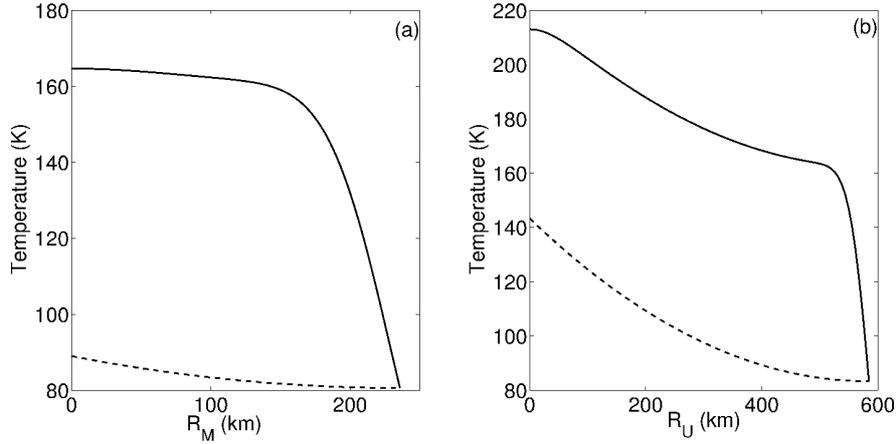}
\caption{Maximal (plain lines) and minimal (dashed lines) profiles for Miranda (a) and Umbriel (b) obtained by the radiogenic heating on 4.6 Gyr.}
\label{figureprofil}
\end{center}
\end{figure}

In the case of Miranda and Umbriel, the results of this simulation is given in Figure \ref{figureprofil} which shows the colder and warmer possible profiles. The scenarios presented 
in the following subsections select the warmer profile as initial temperature profile.

\subsection{Nominal Scenario}

The nominal scenario presented here is a first coupled approach of our problem with realistic thermal parameters and orbital variables.

Let us consider the satellites Miranda and Umbriel with the orbital parameters fixed to the current ones (cf. Table \ref{elementsorbitaux}) except for the semi-major axes considered at the resonance in 
eccentricity $e_M^2$ and the eccentricities fixed to a smaller value. The dynamical evolution inside the resonance is led by the tidal equations (\ref{Kaula}) with
the parameter $(k_2/Q)_p=~ 5.2\ 10^{-5}$ for Uranus \citep{tittemore1988}. For the satellites, the ratio $(k_2/Q)_s$ is directly determined by the thermal code.
The thermal parameters are fixed as resumed in Table \ref{paramthermiques}. The dissipation function $Q$ is computed by the Maxwell model.

\begin{figure}[h!]
  \begin{center}
 \includegraphics[width=140mm]{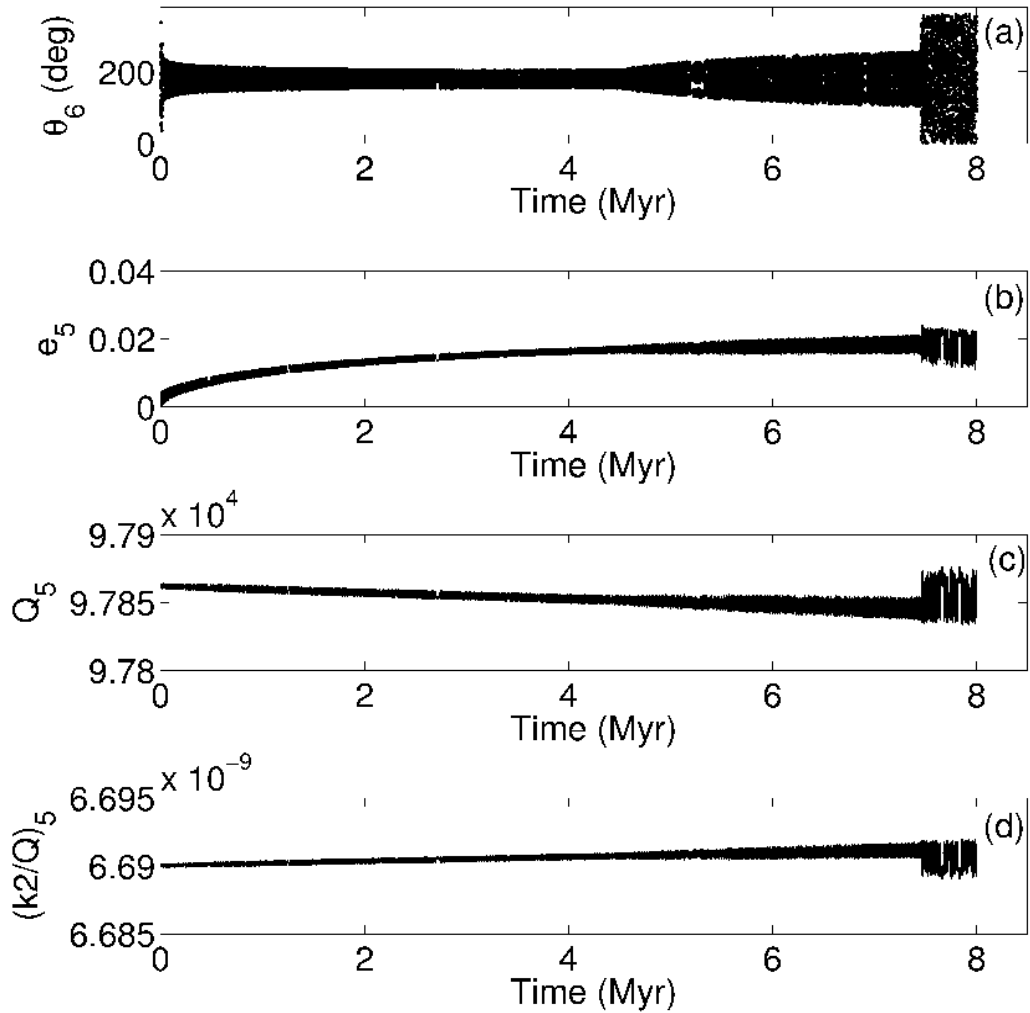}  
  \end{center}
\caption{Results of the coupled orbital-thermal approach with a Maxwell rheology. During the libration of the resonant argument (a), the associated eccentricity
increases (b). Its final value at the exit of the resonance is not maintained and is insufficient to involve a rise of temperature inside Miranda. Consequently the value of $Q$ is rather high with as result a very small value for the ratio $(k_2/Q)_s$, associated with no dissipation inside Miranda. }
\label{T1coupledfigure}
\end{figure}

The results of this first coupled orbital-thermal model are given in Figure \ref{T1coupledfigure}.
The tidal effect on semi-major axes pushes the satellite inside the resonance zone and makes them evolve on a timescale of $6$ Myr.
Inside this resonance $e_M^2$, we observe the libration of the resonant argument $\theta_6$ (a), involving the rise of the eccentricity of Miranda (b). The value of this eccentricity at the exit of the resonance is  
\begin{equation}
 e_5\approx0.02\ .
\end{equation}
This value being quite moderate and, by the tidal synchronization which tends to damp the eccentricities of the orbits, we do not observe any heating 
in Miranda with our choice of initial conditions. With a Maxwell rheology, the viscosity of Miranda is extremely high ($>10^{19}$ Pa s) involving
an important Maxwell characteristic time compared with the orbital periods leading to a full elastic response of the satellite. The rise of eccentricity is not sufficient to dominate the radiogenic heating 
in the heat equation (\ref{heat_equation}) and the dissipation inside the satellite is approximately null. Looking at the tidal equations (\ref{Kaula}), the second term is insignificant
and the tidal evolution is dominated by the dissipation inside the planet: Miranda moves away from the planet and its orbit is circularized.

%

%


An alternative scenario to enhance tidal heating of Miranda would consider another rheological model but, despite smaller values of the dissipation function $Q$ with those models, the viscosity of Miranda stays nevertheless high
to provide significant heating by tides.

It is possible to diminish this viscosity considering a lower melting temperature ($<273$ K) allowing the decrease of this viscosity. Melting temperature as low as $200$ K is possible for Miranda especially if there are  ice clathrates and  elements such as ammonia and salts inside \citep{greenberg1991}.  
Several tests have been performed with different melting temperatures with the three rheological models but none of them gives the heating of Miranda as a result, leading us to conclude that the tidal heating is not effective on Miranda with a classical approach and the chosen parameters.

\begin{figure}[h!]
\begin{center}
\includegraphics[width=140mm]{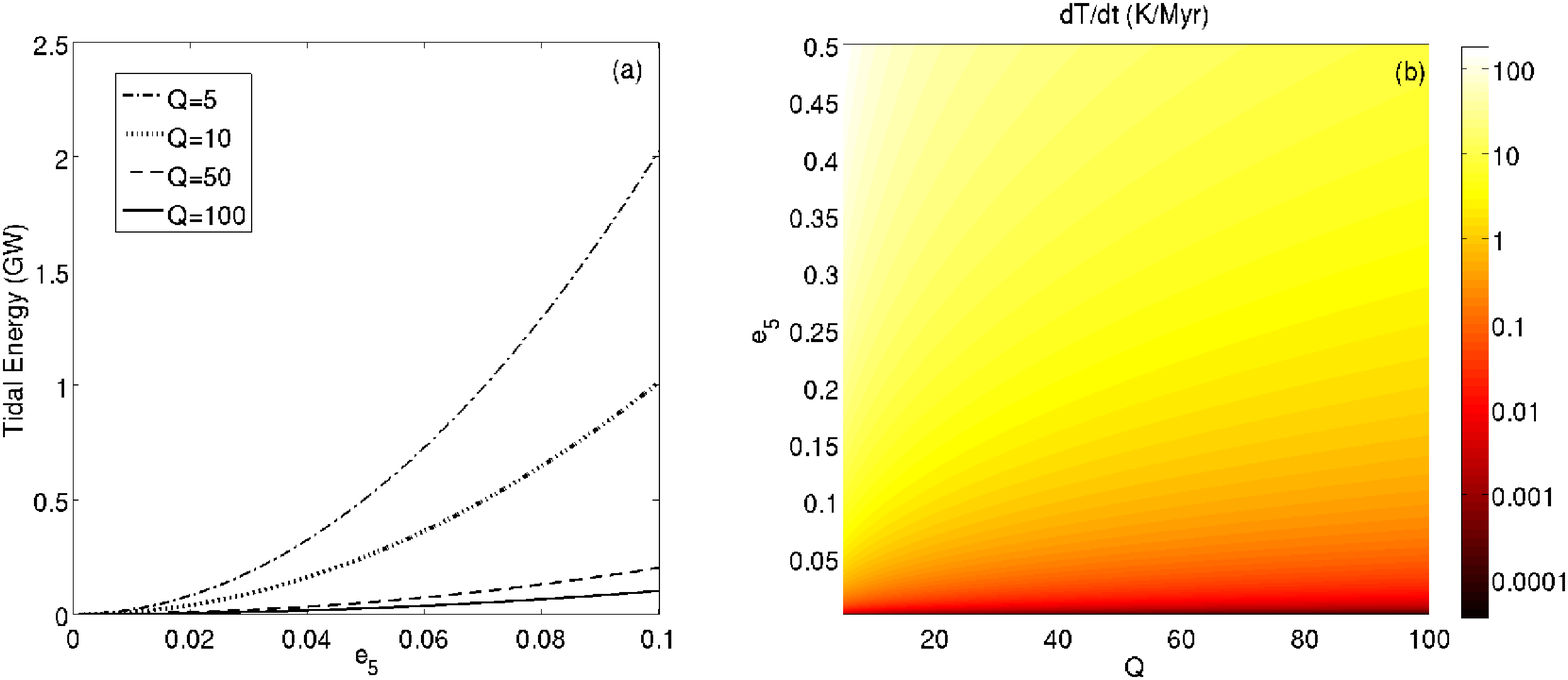}
\end{center}
 \caption{Variation of the tidal energy of dissipation versus the eccentricity for fixed values of $Q$ (a). Variation of the temperature on one million years in a plane eccentricity $e_5$ versus dissipation function $Q$ (b). 
The color scale is logarithmic.}
 \label{ExplicationTest}
 \end{figure}

We confirm this hypothesis by considering the following approximation. Looking at the equation (\ref{energydissipation})
with $f=1$ and $\epsilon=0$, we write : 

\begin{equation}
 C=\frac{21}{2}\ k_2 \frac{G\ m_0^2\ n\ R_s^5}{a^6}\ ,
\end{equation}
where $C$ is constant. The energy by unit of mass evaluated in $e$ and $Q$ is given by :
\begin{equation}
 E= \frac{C}{m_5}\ \frac{e^2}{Q}\ ,
\label{energie}
\end{equation}
where $m_5$ is the mass of Miranda. The total energy is represented in Figure \ref{ExplicationTest} (a). Considering a large dissipation ($Q=5$) we note a value close $1$GW for $e\approx0.06$. Dividing this energy (\ref{energie}) by the specific heat of Miranda, we obtain the variation of internal temperature (K) in one year. 
The figure \ref{ExplicationTest} (b) gives this variation on one million years in a plane eccentricity $e_5$ versus dissipation function $Q$. The color scale is, in this last case, logarithmic.

The approximation given by Figure \ref{ExplicationTest} (b) confirms the impossibility to heat Miranda in the chosen conditions in the first test:
indeed, with the obtained eccentricity lower than $0.05$, the rise of temperature over one million years is less than $1$ K 
for values below $100$ for the dissipation function $Q$. To observe a slight rise of temperature, we need to consider an higher eccentricity close to $0.4$ and 
maintain it on a period of several million years.

We also need quite small values for the dissipation function $Q$, invalidating the use of the Maxwell rheology in our conditions (cf. Figure \ref{dissipationQ}).

\subsection{Extremal Scenario}

The extremal scenario presented in this section considers the rheologies of Burgers and Andrade. We choose
to preserve the physical significance of the thermal parameters to the detriment of plausible dynamical elements in order to obtain a heating
of the satellite interior corresponding to the approximation given in Figure \ref{ExplicationTest} (b).

As the eccentricity of Miranda does not increase sufficiently with the capture and the evolution in the primary resonance $e_M^2$ , we consider a high initial eccentricity ($e_5=0.5$) at the beginning of the simulation
in order to show the effect of the coupling approach. In this case, as explained in \citet{champenois1998}, the eccentricity decreases until an equilibrium value and the exit of the resonance.

\begin{figure}[h!]
  \begin{center}
 \includegraphics[width=140mm]{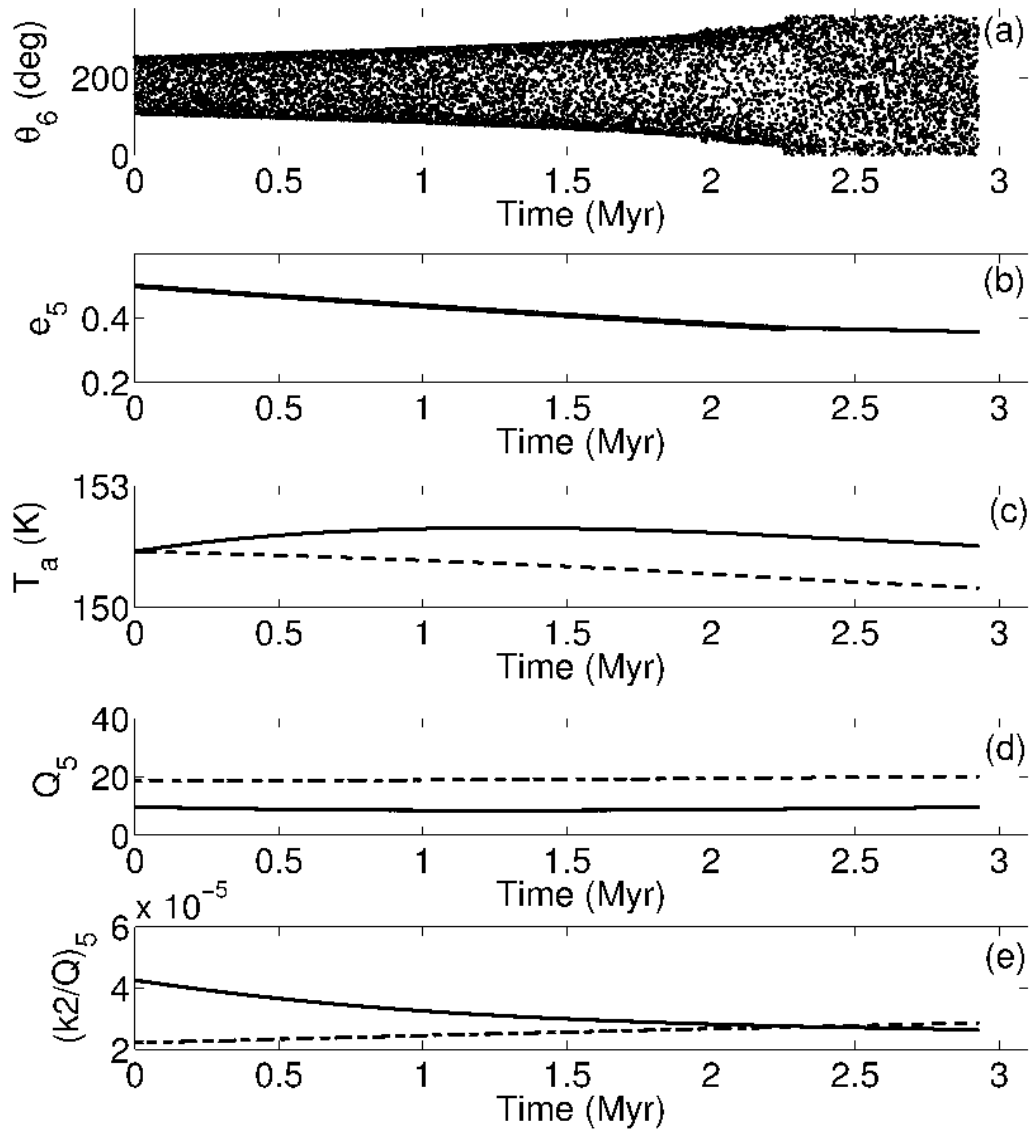}  
  \end{center}
\caption{Results of the coupled orbital-thermal approach with a Burgers (plain lines) and an Andrade (dashed lines) rheologies. During the libration of the resonant argument (a), the associated eccentricity
decreases (b) until an equilibrium value. In both cases, the value of $Q$ is rather small (d) with as a result a reasonable value for the ratio $k_2/Q$ (e), associated with a huge dissipation inside Miranda. The large initial eccentricity,
combined with a melting temperature set to $200$ K, involves a slight increase in the averaged temperature with the Burgers model.}
\label{T2coupledfigure}
\end{figure}

Figure \ref{T2coupledfigure} presents the results with a Burger rheology with $C=50$ (plain lines) and an Andrade model with $\alpha=0.33$ (dashed lines). In both cases, the value of the melting temperature is set to $200$ K to 
diminish the initial viscosity of Miranda and attempt to have a viscous response of the body.
With the two different rheological models, the orbital evolutions are uniform: we observe the exit of the primary resonance illustrated by the end of the libration regime of 
the angle $\theta_6$ (a). The high eccentricity decreases to an equilibrium value close to $0.38$ (b). The coupled thermal evolution differs from the {\it nominal case}
as we observe a strong tidal dissipation for Miranda: the value of $Q$ is rather small and the ratio $(k_2/Q)_s$ is the same order than the ratio $(k_2/Q)_p$. Associated with 
a large eccentricity, the tidal evolution differs from the {\it nominal case} by a predominance of the second term in the equation (\ref{Kaula}): the evolution is now dominated by the 
dissipation inside the satellite. This involves a motion of the satellite towards the planet instead of distancing it and Miranda goes back inside the resonance zone.

The thermal part depends on the chosen model. We observe a slight increase in the averaged temperature with the Burgers model (plain lines (c)). This increase of approximately $1$ K on a million of years corresponds 
to the increase given by the approximation approach illustrated in Figure \ref{ExplicationTest} (b). The value of the dissipation function $Q$ is, in both cases, pretty small (d), involving a reasonable value for the 
ratio $(k_2/Q)_s$ for Miranda (e). However, as the eccentricity value decreases on a quite small scale (the tidal damping prevents the maintaining of the large value for the eccentricity), we observe nevertheless an increase in the dissipation resulting in a cooling of the satellite even with an extremal approach.

%
%
%
%
%

\section{Conclusions and Perspectives}
\label{Conclusions}
In this work, we present on a coupled thermal orbital model based on a 3:1 mean motion resonance in eccentricity 
between Miranda and Umbriel. 

The dynamical module consists of an averaging of the short periodic terms of the complete 3 body problem presented in \citet{verheylewegen2013}.
This approach, using the Hamiltonian formalism, allows us to obtain a 3 dimensional model with the eccentricities and the inclinations
of the orbits. It is coupled to a thermal module via the tidal effect.

The thermal module presents the resolution of the heat equation for a one dimensional homogenous sphere composed by a mixture of silicates and ices.
We considered a temperature dependent viscosity and different rheological models. The main mode of heat transfer is conduction since in the simulations, 
the internal heat is not sufficient to start convection. The temperature inside the satellite can increase due to the decay of radiogenic elements or with the tidal dissipation depending on orbital elements. 
The thermal parameters depend on the averaged temperature inside the satellites.

A coupled approach of the orbital and internal dynamics has never been applied in the case of the main satellites of Uranus. The solutions depend
essentially on the chosen rheological model and on the eccentricities of the satellites. Regardless the eccentricity, the Maxwell model is ineffective with our 
conductive approach given the important values of the dissipation function $Q$ for the satellites. The other two models (Burgers and Andrade) give 
values of $Q$ around $100$ for Miranda. However its final eccentricity $e\approx 0.02$, obtained by the capture into the resonance 3:1 with Umbriel, is not maintained due to 
the tidal damping acting on a short time-scale and Miranda cools down. 

An alternative scenario is studied with higher initial eccentricity leading to a slight heating of Miranda on one million of years. Orbital evolution in this extreme case results in 
a decrease of the eccentricity until an equilibrium value associated with a diminution of the dissipated tidal energy is reached after which Miranda cools down again.

These results of cooling Miranda do not exclusively depend on the eccentricities and on the rheological models but are also associated with our model's choice of a 
conductive homogenous sphere, initially at a low temperature. They do not invalidate any possibility of internal heating for Miranda.

By comparison, Enceladus and Mimas on the system of Saturn have many similarities with the Uranian satellites and, in particular, with Miranda. Mimas, despite its large eccentricity and 
closer distance to Saturn, does not show past or present geological activity whereas Enceladus is geologically very active today. 

Indeed, Cassini's composite infrared spectrometer of Enceladus' south polar terrain, which is marked by linear fissures, indicates that the internal heat-generated power is about $15.8$ GW \citep{howett2011}
or, more recently, about $4.7$ GW \citep{spencer2013}. 
Water-rich plume venting from the moon's south polar region associated with the large internal heating makes very likely the presence of liquid water below the Enceladus surface. 

A south polar sea between the moon's outer ice shell and its rocky interior would increase the efficiency of the tidal heating by allowing greater tidal distortions of the 
ice shell. The difference between the recent geological history of Mimas, Miranda and Enceladus is associated with their temperature dependent material properties  such as viscosity  
$\eta$ and dissipation factor $Q$. An increase of  internal temperatures due to for instance radiogenic heating, orbital resonance or catastrophic events would enhance tidal heating since 
$\eta$ and $Q$ would decrease exponentially with increasing temperatures.

Despite Enceladus being further from Saturn and lower eccentricity, its current high-energy thermal state is likely linked to its thermal evolution. The high internal heating is generally attributed  to the tidal heating 
enhanced by orbital resonances. The heating in Enceladus in an equilibrium resonant configuration with other Saturnian satellites is studied by \citet{meyer2007}. They showed that equilibrium tidal heating cannot 
account for the heat that is observed to be coming from Enceladus unless we consider high dissipation in Saturn \citep{lainey2012}. The equilibrium heating in possible past resonances likewise cannot explain prior resurfacing events. 
While the exact source and mechanism of Enceladus' internal heating is currently not known, it provides a good analogy for Uranian satellites that we will consider in future studies.  

\section*{Acknowledgments}
The work of Emilie Verheylewegen is supported by an FNRS PhD Fellowship. The work of Beno\^it Noyelles is supported by an FNRS Postdoctoral Research Fellowship. \"Ozg\"ur Karatekin 
benefits from the financial support of the Belgian PRODEX, managed by the ESA, in collaboration with the Belgian Federal Science Policy Office.
This research used resources of the "Plateforme Technologique de Calcul Intensif (PTCI)" (http://www.ptci.unamur.be) located at the University of Namur, Belgium, which is supported 
by the F.R.S.-FNRS. The authors want to thank Attilio Rivoldini for his valuable advice and discussion for the thermal module and Andr\'e F\"uzfa for his advice in EDP resolution.


%
%





\bibliographystyle{apalike}
\bibliography{biblio}

\begin{thebibliography}{}

\bibitem[Andrade, 1910]{andrade1910}
Andrade, E. (1910).
\newblock On the viscous flow in metals, and allied phenomena.
\newblock In {\em Royal Society Proceedings}, volume~84, pages 1--12. The Royal
  Society of London.

\bibitem[Archinal et~al., 2011]{Archinal2009}
Archinal, B., A'Hearn, M., Bowell, E., et~al. (2011).
\newblock Report of the {IAU} {W}orking group on {C}artographic {C}oordinates
  and {R}otational {E}lements: 2009.
\newblock {\em Celestial Mechanics and Dynamical Astronomy}, 109:101--135.

\bibitem[{Brown} and {Clark}, 1984]{brown1984}
{Brown}, R. and {Clark}, R. (1984).
\newblock {Surface of Miranda - Identification of water ice}.
\newblock {\em Icarus}, 58:288--292.

\bibitem[Brown et~al., 1991]{brown1991}
Brown, R., Johnson, T., Synnott, S., Anderson, J., and Jacobson, R. (1991).
\newblock {P}hysical properties of the {U}ranian {S}atellites.
\newblock In {\em {U}ranus}, {Space Science Series}, pages 513--527. The
  University of Arizona Press.

\bibitem[Castillo-Rogez et~al., 2007]{castillo2007}
Castillo-Rogez, J., Matson, D., Sotin, S., Johnson, T., Lunine, J., and Thomas,
  P. (2007).
\newblock Iapetus' geophysics: {R}otation rate, shape, and equatorial ridge.
\newblock {\em Icarus}, 190:179--202.

\bibitem[Champenois, 1998]{champenois1998}
Champenois, S. (1998).
\newblock {\em {D}ynamique de la r{\'e}sonance entre {M}imas et {T}{\'e}thys,
  premier et troisi{\`e}me satellites de {S}aturne}.
\newblock PhD thesis, Observatoire de Paris, France.

\bibitem[Cincotta and Sim{\`o}, 2000]{cincotta2000}
Cincotta, P. and Sim{\`o}, C. (2000).
\newblock Simple tools to study global dynamics in non-axisymmetric galactic
  potentials - {I}.
\newblock {\em Astronomy and Astrophysics}, 147:205--228.

\bibitem[Dermott et~al., 1988]{dermott1988}
Dermott, S., Malhotra, R., and Murray, C. (1988).
\newblock {D}ynamics of the {U}ranian and {S}aturnian satellite systems- {A}
  chaotic route to melting {M}iranda?
\newblock {\em Icarus}, 76:295--334.

\bibitem[Douce, 2011]{douce2011}
Douce, E. (2011).
\newblock {\em {T}hermodynamics of the {E}arth and {P}lanets}.
\newblock Cambridge University Press, New York, 1st edition.

\bibitem[Duriez, 1977]{duriez1977}
Duriez, L. (1977).
\newblock Th{\'e}orie {G}{\'e}n{\'e}rale {P}lan{\'e}taire en {V}ariables
  {E}lliptiques. {I}. {D}{\'e}veloppement des {E}quations.
\newblock {\em Astronomy and Astrophysics}, 54:93--112.

\bibitem[{Efroimsky}, 2012]{efroimsky2012}
{Efroimsky}, M. (2012).
\newblock {Tidal Dissipation Compared to Seismic Dissipation: In Small Bodies,
  Earths, and Super-Earths}.
\newblock {\em The Astrophysical Journal}, 746:150--170.

\bibitem[Greenberg et~al., 1991]{greenberg1991}
Greenberg, S., Croft, J., Janes, D., Kargel, J., Lebofsky, L., Lunine, J.,
  Marcialis, H., Melosh, G., G.W., O., and Strom, R. (1991).
\newblock {M}iranda.
\newblock In Bergstralh, J., Miner, E., and Matthews, S., editors, {\em
  {U}ranus}. The University of Arizona Press.

\bibitem[{Howett} et~al., 2011]{howett2011}
{Howett}, C.~J.~A., {Spencer}, J.~R., {Pearl}, J., and {Segura}, M. (2011).
\newblock {High heat flow from Enceladus' south polar region measured using
  10-600 cm$^{-1}$ Cassini/CIRS data}.
\newblock {\em Journal of Geophysical Research (Planets)}, 116:3003.

\bibitem[Hussmann et~al., 2006]{hussmann2006}
Hussmann, H., Sohl, F., and Spohn, T. (2006).
\newblock {S}ubsurface oceans and deep interiors of medium-sized outer planet
  satellites and large trans-neptunian objects.
\newblock {\em Icarus}, 185:258--273.

\bibitem[Hussmann et~al., 2009]{hussmann2009}
Hussmann, H., Sotin, C., and Lunine, J. (2009).
\newblock {\em Interiors and {E}volution of {I}cy {S}atellites}, volume~10 of
  {\em Planets and Moons: Treatrise on Geophysics}, chapter~15.
\newblock Elsevier, Germany.

\bibitem[{Jacobson}, 2007]{jacobson2007}
{Jacobson}, R. (2007).
\newblock {The Gravity Field of the Uranian System and the Orbits of the
  Uranian Satellites and Rings}.
\newblock In {\em AAS/Division for Planetary Sciences Meeting Abstracts \#39},
  volume~39, page 453. Bulletin of the American Astronomical Society.

\bibitem[Karato, 1998]{karato1998}
Karato, S. (1998).
\newblock {\em {Deformation of Earth Materials: An Introduction to the Rheology
  of Solid Earth}}.
\newblock Cambridge University Press, New York.

\bibitem[{Kargel} and {Lewis}, 1993]{kargel1993}
{Kargel}, J. and {Lewis}, J. (1993).
\newblock {The Composition and Early Evolution of Earth}.
\newblock {\em Icarus}, 105:1--25.

\bibitem[Lainey et~al., 2012]{lainey2012}
Lainey, V., Karatekin, {\"O}., Desmars, J., et~al. (2012).
\newblock {S}trong tidal dissipation in {S}aturn and constraints on
  {E}nceladus' thermal state from astrometry.
\newblock {\em The Astrophysical Journal}, 752:14--23.

\bibitem[Laskar, 1993]{laskar1993}
Laskar, J. (1993).
\newblock Frequency analysis of a dynamical system.
\newblock {\em Celestial Mechanics and Dynamical Astronomy}, 56:191--196.

\bibitem[Laskar and Jacobson, 1987]{laskar1987}
Laskar, J. and Jacobson, R. (1987).
\newblock {GUST}86. {A}n analytical ephemeris of the {U}ranian satellites.
\newblock {\em Astronomy and Astrophysics}, 188:212--224.

\bibitem[Malhotra, 1991]{malhotra1991}
Malhotra, R. (1991).
\newblock {T}idal origin of the {L}aplace {R}esonance and the resurfacing of
  {G}anymede.
\newblock {\em Icarus}, 94:399--412.

\bibitem[Malhotra and Dermott, 1990]{malhotradermott1990}
Malhotra, R. and Dermott, S. (1990).
\newblock {T}he role of secondary resonances in the orbital history of
  {M}iranda.
\newblock {\em Icarus}, 85:444--480.

\bibitem[Mavko et~al., 2009]{mavko2009}
Mavko, G., Mukerji, T., and Dvorkin, J. (2009).
\newblock {\em {The Rock Physics Handbook : Tools for Seismic Analysis of
  Porous Media}}.
\newblock Cambridge University Press, New York, 2{\`e}me edition.

\bibitem[{Meyer} and {Wisdom}, 2007]{meyer2007}
{Meyer}, J. and {Wisdom}, J. (2007).
\newblock {Tidal heating in Enceladus}.
\newblock {\em Icarus}, 188:535--539.

\bibitem[Murray and Dermott, 1999]{murray1999}
Murray, C. and Dermott, S. (1999).
\newblock {\em {S}olar {S}ystem {D}ynamics}.
\newblock Cambridge University Press, New York.

\bibitem[{Noyelles}, 2010]{noyelles2010}
{Noyelles}, B. (2010).
\newblock {Theory of the rotation of Janus and Epimetheus}.
\newblock {\em Icarus}, 207:887--902.

\bibitem[Ojakangas and Stevenson, 1986]{ojakangas1986}
Ojakangas, G. and Stevenson, D. (1986).
\newblock {E}pisodic {V}olcanism of {T}idally {H}eated {S}atellites with
  {A}pplication to {I}o.
\newblock {\em Icarus}, 66:341--358.

\bibitem[Pappalardo et~al., 1997]{pappalardo1997}
Pappalardo, R., Reynolds, S., and Greeley, R. (1997).
\newblock {E}xtensional tilt blocks on {M}iranda: {E}vidence for an upwelling
  origin of {A}rden {C}orona.
\newblock {\em Journal of Geophysical Research}, 102:13369--13379.

\bibitem[Parmentier and Zuber, 2007]{parmentier2007}
Parmentier, E. and Zuber, M. (2007).
\newblock {E}arly evolution of {M}ars with mantle compositional stratification
  or hydrothermal crustal cooling.
\newblock {\em Journal of Geophysical research}, 112:2007--2018.

\bibitem[{Peale}, 1988]{peale1988b}
{Peale}, S. (1988).
\newblock {Speculative histories of the Uranian satellite system}.
\newblock {\em Icarus}, 74:153--171.

\bibitem[{Peale}, 1999]{peale1999}
{Peale}, S. (1999).
\newblock {Origin and Evolution of the Natural Satellites}.
\newblock {\em Annual Review of Astronomy and Astrophysics}, 37:533--602.

\bibitem[{Peltier}, 1974]{peltier1974}
{Peltier}, W. (1974).
\newblock {The impulse response of a Maxwell earth.}
\newblock {\em Reviews of Geophysics and Space Physics}, 12:649--669.

\bibitem[{Plescia}, 1987]{plescia1987}
{Plescia}, J. (1987).
\newblock {Geological terrains and crater frequencies on Ariel}.
\newblock {\em Nature}, 327:201--204.

\bibitem[{Plescia}, 1988]{plescia1988}
{Plescia}, J. (1988).
\newblock {Cratering history of Miranda - Implications for geologic processes}.
\newblock {\em Icarus}, 73:442--461.

\bibitem[{Rambaux} et~al., 2010]{rambaux2010}
{Rambaux}, N., {Castillo-Rogez}, J., {Williams}, J., and {Karatekin}, {\"O}.
  (2010).
\newblock {Librational response of Enceladus}.
\newblock {\em Geophysical Research Letters}, 37:4202--4207.

\bibitem[Reeh et~al., 2003]{reeh2003}
Reeh, N., Christensen, E., Meyer, C., and Olesen, O. (2003).
\newblock {Tidal bending of glaciers: A linear viscoelastic approach}.
\newblock {\em Annals of Glaciology}, 37:83--89.

\bibitem[{Robuchon} et~al., 2010]{robuchon2010}
{Robuchon}, G., {Choblet}, G., {Tobie}, G., {\v{C}adek}, O., {Sotin}, C., and
  {Grasset}, O. (2010).
\newblock {Coupling of thermal evolution and despinning of early Iapetus}.
\newblock {\em Icarus}, 207:959--971.

\bibitem[{Schubert} et~al., 2010]{schubert2010}
{Schubert}, G., {Hussmann}, H., {Lainey}, V., et~al. (2010).
\newblock {Evolution of Icy Satellites}.
\newblock {\em Space Science review}, 153:447--484.

\bibitem[Schubert et~al., 2001]{schubert2001}
Schubert, G., Turcotte, D., and Olson, P. (2001).
\newblock {\em {M}antle convection in the {E}arth and {P}lanets}.
\newblock Cambridge University Press, New York, 1st edition.

\bibitem[{Shoji} et~al., 2013]{shoji2013}
{Shoji}, D., {Hussmann}, H., {Kurita}, K., and {Sohl}, F. (2013).
\newblock {Ice rheology and tidal heating of Enceladus}.
\newblock {\em Icarus}, 226:10--19.

\bibitem[{Showman} et~al., 1997]{showman1997}
{Showman}, A., {Stevenson}, D., and {Malhotra}, R. (1997).
\newblock {Coupled Orbital and Thermal Evolution of Ganymede}.
\newblock {\em Icarus}, 129:367--383.

\bibitem[{Smith} et~al., 1986]{smith1986}
{Smith}, B., {Soderblom}, L., {Beebe}, R., et~al. (1986).
\newblock {Voyager 2 in the Uranian system - Imaging science results}.
\newblock {\em Science}, 233:43--64.

\bibitem[{Spencer} et~al., 2013]{spencer2013}
{Spencer}, J.~R., {Howett}, C.~J., {Verbiscer}, A.~J., {Hurford}, T.~A.,
  {Segura}, M., and {Spencer}, D.~C. (2013).
\newblock {A New Estimate of the Power Emitted by Enceladus' Tiger Stripes}.
\newblock In {\em AAS/Division for Planetary Sciences Meeting Abstracts},
  volume~45 of {\em AAS/Division for Planetary Sciences Meeting Abstracts},
  page \#403.03.

\bibitem[{Strobell} and {Masursky}, 1987]{strobell1987}
{Strobell}, M.~E. and {Masursky}, H. (1987).
\newblock {New Features Named on the Moon and Uranian Satellites}.
\newblock In {\em Lunar and Planetary Institute Science Conference Abstracts},
  volume~18, page 964. Lunar and Planetary Inst. Technical Report.

\bibitem[Thomas, 1988]{thomas1988}
Thomas, P. (1988).
\newblock {R}adii, shapes, and topography of the satellites of {U}ranus from
  limb coordinates.
\newblock {\em Icarus}, 73:427--441.

\bibitem[Tittemore and Wisdom, 1988]{tittemore1988}
Tittemore, W. and Wisdom, J. (1988).
\newblock {T}idal {E}volution of the {U}ranian {S}atellites- {I}. {P}assage of
  {A}riel and {U}mbriel through the $5:3$ {M}ean-{M}otion {C}ommensurability.
\newblock {\em Icarus}, 74:172--230.

\bibitem[Tittemore and Wisdom, 1989]{tittemore1989}
Tittemore, W. and Wisdom, J. (1989).
\newblock {T}idal {E}volution of the {U}ranian {S}atellites- {II}. {A}n
  {E}xplanation of the {A}nomalously {H}igh {O}rbital {I}nclination of
  {M}iranda.
\newblock {\em Icarus}, 78:63--89.

\bibitem[Tittemore and Wisdom, 1990]{tittemore1990}
Tittemore, W. and Wisdom, J. (1990).
\newblock {T}idal {E}volution of the {U}ranian {S}atellites- {III}. {E}volution
  through the {M}iranda-{U}mbriel $3:1$, {M}iranda-{A}riel $5:3$, and
  {A}riel-{U}mbriel $2:1$ {M}ean-{M}otion {C}ommensurabilities.
\newblock {\em Icarus}, 85:394--443.

\bibitem[Verheylewegen et~al., 2013]{verheylewegen2013}
Verheylewegen, E., Noyelles, B., and Lema{\^i}tre, A. (2013).
\newblock {A} numerical exploration of {M}iranda's dynamical history.
\newblock {\em The Monthly Notices of the Royal Astronomical Society},
  435:1776--1787.

\bibitem[Yoder and Peale, 1981]{Yoder1981}
Yoder, C. and Peale, S. (1981).
\newblock {T}he tides of {I}o.
\newblock {\em Icarus}, 47:1--35.

\end{thebibliography}







\end{document}